\documentclass[a4paper, 11pt
]{amsart}
\usepackage[latin1]{inputenc}
\usepackage[english]{babel}
\usepackage[T1]{fontenc}
\usepackage{amsmath,amssymb, amsfonts}
\usepackage{amsthm}
\usepackage{latexsym}

\usepackage{epsfig}
\usepackage{graphicx}
\usepackage{color}

\usepackage{syntonly}

\usepackage[all, curve, frame, arc, color,
 dvips 
]{xy} \xyoption{tile}

\usepackage[dvips]{hyperref}

\renewcommand{\cong}{\simeq}

\newcommand{\ph}[0]{\varphi}

\renewcommand{\i}{{\mathrm{i}}}

\newcommand{\br}[1]{G_{#1}}

\newcommand{\R}{{\mathbb{R}}}
\newcommand{\C}{{\mathbb{C}}}
\newcommand{\Z}{{\mathbb{Z}}}
\newcommand{\Q}{{\mathbb{Q}}}
\newcommand{\N}{{\mathbb{N}}}

\newcommand{{\abin}}[2]{  \begin{array}{c} #1 \\ #2
\end{array}  }
\newcommand{\oC}[0]{L^1}
\newcommand{\cF}[0]{\mathcal{F}}

\newcommand{\di}{{\mathrm{d}}}

\newcommand{\inv}{I}
\newcommand{\anti}{K}
\newcommand{\mb}[1]{\mathbf{#1}}

\newcommand{\card}[1]{|#1|}

\newcommand{\qbin}[2]{ \left[ \begin{array}{c} #1 \\ #2
    \end{array} \right] }

\renewcommand{\mod}[1]{\,\mathrm{mod}(#1)}

\newenvironment{dm}[1][Proof.]{\begin{trivlist}
\item[\hskip
    \labelsep {\bfseries #1}]}{\end{trivlist}}

\newcommand{\tss}[0]{\supset}
\newcommand{\sst}[0]{\subset}
\newcommand{\pmu}[0]{{\pm 1}}

\newcommand{\vs}[0]{\vspace{3mm}}

\renewcommand{\t}[1]{\widetilde{#1}}

\newcommand{\too}{\longrightarrow}

\def\qed{\ifmmode   \textrm{ \qquad \qquad $\Box$} \else{\unskip\nobreak\hfil
\penalty50\hskip1em\null\nobreak\hfil $\Box$
\parfillskip=0pt\finalhyphendemerits=0\endgraf}\fi}
\newtheorem{df}{Definition}[section]

\newtheorem{prop}[df]{Proposition}

\newtheorem{lemma}[df]{Lemma}
\newtheorem{theorem}[df]{Theorem}

\begin{document}
\title[The $K(\pi, 1)$ problem for $\widetilde{B}_n$ and its cohomology]{The $K(\pi, 1)$ problem for the affine Artin group of type
  $\t{B}_n$ and its cohomology}
\date{\today}

\author[F. Callegaro]{Filippo Callegaro}
\address{Scuola Normale Superiore\\ P.za dei Cavalieri, 7, Pisa, Italy}
\email{f.callegaro@sns.it}

\author[D. Moroni]{Davide Moroni}
\address{Dipartimento di Matematica ``G.Castelnuovo''\\
P.za A.~Moro, 2, Roma, Italy -and- ISTI-CNR\\ Via G.~Moruzzi, 3,  Pisa, Italy}
\email{davide.moroni@isti.cnr.it}

\author[M. Salvetti]{Mario Salvetti}
\address{Dipartimento di Matematica ``L.Tonelli''\\ Largo B.~Pontecorvo, 5, Pisa, Italy}
\email{salvetti@dm.unipi.it}
\thanks{The third author is partially supported by M.U.R.S.T. 40\%}


\begin{abstract}

In this paper we prove that the complement to the affine complex 
arrangement of type $\widetilde{B}_n$ is a $K(\pi, 1)$ space.  
We also compute the cohomology of the affine Artin group 
$G_{\widetilde{B}_n}$ (of type $\widetilde{B}_n$)  
with coefficients over several interesting local systems. In particular,
we consider the module $\Q[q^{\pm 1}, t^{\pm 1}],$ where the
first $n$-standard generators of $G_{\widetilde{B}_n}$ act  
by $(-q)$-multiplication while the last generator acts by
$(-t)$-multiplication. Such representation  
generalizes the analog $1$-parameter representation related 
to the bundle structure over the complement to the discriminant
hypersurface,  endowed  with the monodromy action of the  
associated Milnor fibre. The cohomology of  
$G_{\widetilde{B}_n}$ with trivial coefficients is derived from the
previous one.

\end{abstract}


\maketitle 



\section{Introduction}

Let $(W,S)$ be a Coxeter system, so a presentation for $W$ is
$$<s\in S\ |\ (ss')^{m(s,s')}=1>$$
where $m(s,s')\in \mathbb N_{\geq 2}\cup \{\infty\}$ for $s\neq s'$ and
  $m(s,s)=1$ (see \cite{bourbaki},
\cite{humphreys}).

The Artin group $\br{W}$ associated to $(W,S)$ is
the extension of $W$ given by the presentation (see \cite{Br-Sa})
$$<g_s,\ s\in S\ |\ g_sg_{s'}g_s...=g_{s'}g_sg_{s'}... \ (s\neq s',\
m(s,s')\ \mbox{factors})>.$$ One says that an Artin group $\br{W}$ is of
\emph{finite type} when $W$ is finite. We are interested in
\emph{finitely generated} Artin groups, that is when $S$ is finite.
In this case, $W$ can be geometrically represented as a linear
reflection group in $\R^n$ (for example, by using the \emph{Tits
representation} of $W,$ see \cite{bourbaki}). Let $\mathcal{A}^\R$
be the arrangement of hyperplanes given by the mirrors of the
reflections in $W$ and let its complement be
$\mb{Y}(\mathcal{A}^\R):=\R^n \setminus \bigcup_{\mb{H}^\R \in \mathcal{A}^\R}
\mb{H}^\R.$ The connected components of the complement
${\mb{Y}}(\mathcal{A^\R})$ are called the \emph{chambers} of
$\mathcal{A^\R}$.

Consider (for finite type) the arrangement $\mathcal{A}$ in $\C^n$
obtained by complexifying the hyperplanes of $\mathcal{A^\R}$ and
let $\mb{Y}(\mathcal{A})$ be its complement. We have an induced action of
$W$ on $\mb{Y}(\mathcal{A})$ and it turns out that the  {\it orbit space}
$\mb{Y}(\mathcal{A})/W$ has the Artin group $\br{W}$ as fundamental group
(see \cite{brieskorn}). Moreover, it follows
from a Theorem by Deligne (\cite{deligne_imm}) that
$\mb{Y}(\mathcal{A})/W$ is a $K(\pi,1)$ space. Indeed the Theorem
concerns a more general situation. Recall that a real arrangement
$\mathcal{A}^\R$ is said to be \emph{simplicial} if all its chambers
consist of simplicial cones; reflection arrangements are known to be
simplicial \cite{bourbaki}.

\begin{theorem}\label{thm:deligne}\cite{deligne_imm}
Let $\mathcal{A}^\R$ be a finite central arrangement and let
$\mb{Y}(\mathcal{A})$ be the complement of its complexification. If
$\mathcal{A^\R}$ is simplicial, then $\mb{Y}(\mathcal{A})$ is a
$K(\pi,1)$ space. \qed
\end{theorem}

Infinite type Artin groups are represented (by Tits representation;
see also \cite{vinberg} for more general constructions) as groups of
linear, not necessarily orthogonal, reflections w.r.t. the walls of a
polyhedral cone $C$ of maximal dimension in $\mb{V}=\R^n$. It can be
shown that the union $U=\bigcup_{w \in W}wC$ of $W$-translates of
$C$ is a convex cone and that $W$ acts properly on the interior
$U^{0}$ of $U$. We may now rephrase the construction used in the
finite case as follows. Let $\mathcal{A}$  be the complexified
arrangement of the mirrors of the reflections in $W$ and consider
$I:=\{v \in \mb{V} \otimes \C \,|\, \Re(v)\in U^0\}$. Then $W$ acts
freely on $\mb{Y}=I\setminus \bigcup_{\mb{H} \in \mathcal{A}}\mb{H}$ and we can
form the orbit space  $\mb{X}:=\mb{Y}/W$. It is known (\cite{vanderlek}; see
also \cite{salvettiArtin}) that $\br{W}$ is indeed the fundamental group of
$\mb{X}$, but in general it is only conjectured that $\mb{X}$ is a $K(\pi,1)$.
This conjecture is known to be true for: 1) Artin groups of large
type (\cite{hendriks}), 2) Artin groups satisfying the FC condition
(\cite{charney_davis}) and 3) for the affine Artin group of type
$\t{A}_n, \t{C}_n$ (\cite{okonek}). In this note, we extend this
result to the affine Artin group of type $\t{B}_n$, showing:

\begin{theorem} \label{thm:main}
$\mb{Y}(\t{B}_n)$ and, hence, $\mb{X}(\t{B}_n)$ are  $K(\pi, 1)$ spaces.
\end{theorem}

The idea of proof can be described in few words: up to a $\C^*$
factor, the orbit space is presented (through the exponential map)
as a covering of the complement to a finite simplicial arrangement,
so we apply Theorem \ref{thm:deligne}.

We just digress a bit on the peculiarity of affine Artin groups.
In this case the associated Coxeter group is an affine Weyl group
$W_a$ and, as such, it can be geometrically represented as a group
generated by affine (orthogonal) reflections in a real vector space.
This geometric representation and that given by the Tits cone are
linked in a precise manner; indeed it turns out that $U_0$ for an
affine Weyl group is an open half space in $\mb{V}$ and that $W_a$ acts
as a group of affine orthogonal reflections on a hyperplane section
$E$ of $U_0$. The representation on $E$ coincides with the geometric
representation and $\mb{Y}(W_a)$ is homotopic to the
complement of the complexified affine reflection arrangement.

Our second main result is the computation of the cohomology of the
group $\br{\t{B}_n}$ (so, by Theorem \ref{thm:main}), of
$\mb{X}(\t{B}_n)$) with local coefficients. We consider the
$2$-parameters representation of $\br{\t{B}_n}$ over the ring
$\Q[q^{\pm 1},t^{\pm1}]$ and over the module $\Q[[q^\pmu, t^\pmu]]$ defined by sending the standard generator corresponding to the last node of the Dynkin diagram to
$(-t)-$multiplication and the other standard generators to
$(-q)-$multiplication (minus sign is only for technical reasons).
Such representations are quite natural to be considered: they
generalize the analog $1$-parameter representations that (for finite
type) correspond to considering the structure of bundle over the
complement of the discriminant hypersurface in the orbit space and
the monodromy action on the cohomology of the associated
Milnor fibre (see for example \cite{Fr}, \cite{cs}). We explain in
Section \ref{par:shift} various relations between these cohomologies
and the cohomology of the commutator subgroup of $\br{\t{B}_n}.$

The main tool to perform computations is an algebraic complex which
was discovered in \cite{salvettiArtin},
\cite{salvettiDeconciniNotaArtin}
 by using topological
methods (and independently, by algebraic methods in \cite{Sq}).
The cohomology factorizes into two parts (see also
\cite{arithmetic_artin}) : the {\it invariant} part reduces to that of
the Artin group of finite type $B_n,$ whose $2$-parameters
cohomology was computed in \cite{cohom_B_n}; for the {\it
anti-invariant} part we use suitable filtrations and the associated
spectral sequences.

Let $\ph_d$ be the $d$-th cyclotomic polynomial in the variable $q$. We define the quotient rings
$$\{ 1 \}_
i= \Q[q^\pmu, t^\pmu]/(1+tq^i)$$
$$\{ d \}_i = \Q[q^\pmu, t^\pmu]/(\ph_d, 1+tq^i)$$
$$\{ \{ d \} \}_j = \Q[q^\pmu, t^\pmu]/(\ph_d, \prod_{i=o}^{d-1} 1+tq^i)^j.$$
The final result is the following one:

%
\begin{theorem}\label{thm:main2}
The cohomology $H^{n-s}({G}_{\widetilde{B}_n},
\Q[[q^{\pmu},t^{\pmu}]])$ is given by
$$
\begin{array}{crl}
\Q[[q^\pmu, t^\pmu]] &
\mbox{for} \quad  s=&0 \\
& \\
\displaystyle \bigoplus_{h>0} \{ \{ 2h \} \}_{f(n,h)} &
\mbox{for} \quad  s=&1 \\
\end{array}
$$
\begin{align*}
\displaystyle \bigoplus_{\begin{array}{c}\scriptstyle h>2 \\
\scriptstyle i \in I(n,h) \end{array}} \{2h\}_i^{c(n,h,s)} &\oplus
\bigoplus_{\begin{array}{c} \scriptstyle d \mid n \\\scriptstyle 0
\leq i \leq d-2
\end{array}
} \{ d \}_i \oplus \{ 1 \}_{n-1} &
\mbox{for } \quad s =& 2\\
\displaystyle \bigoplus_{\begin{array}{c} \scriptstyle h>2 \\
\scriptstyle i \in I(n,h) \end{array}} \{2h\}_i^{c(n,h,s)} &\oplus
\bigoplus_{
\begin{array}{c}
\scriptstyle d \mid n \\ \scriptstyle 0 \leq i \leq d-2 \\
\scriptstyle d \leq \frac{n}{j+1}
\end{array} }
\{ d \}_i &
\mbox{ for } \quad s =& 2+2j \\
\displaystyle \bigoplus_{\begin{array}{c} \scriptstyle h>2 \\
\scriptstyle i \in I(n,h) \end{array}} \{2h\}_i^{c(n,h,s)} &\oplus
\bigoplus_{
\begin{array}{c}
\scriptstyle d \nmid n \\
\scriptstyle d \leq \frac{n}{j+1}
\end{array}
} \{ d \}_{n-1} & \mbox{for } \quad s =& 3+ 2j
\end{align*}
where $c(n,h,s) = max(0,\lfloor \frac{n}{2h}\rfloor-s)$, $ f(n,h)= \lfloor \frac{n+h-1}{2h} \rfloor$ and
$I(n,h) = \{n, \ldots, n+h-2 \}$ if $n \equiv 0, 1, \ldots, h
\mod{2h}$ and $I(n,h) = \{n+h-1, \ldots, n+2h-1 \} $ if $n \equiv
h+1, h+2, \ldots, 2h-1 \mod{2h}$.
\end{theorem}

As a corollary we also derive the cohomology with trivial
coefficients of $\br{\t{B}_n}$ (Theorem \ref{thm:main3})

The paper is organized as follows. In Section \ref{sec:pre} we
recall some result and notations about Coxeter and Artin groups,
including a $2$-parameters Poincar\'e series which we need in the
boundary operators of the above mentioned algebraic complex. In
Section \ref{sec:kpi_1} we prove Theorem \ref{thm:main}. In
Section \ref{sec:cohom} we use a suitable filtration of the
algebraic complex, reducing computation of the cohomology mainly
to:
\begin{itemize}
    \item calculation of generators of certain subcomplexes for the
Artin group of type  $D_n$ (whose cohomology was known from
\cite{arithmetic_artin}, but we need explicit suitable
generators);
    \item analysis of the associated spectral sequence to
deduce the cohomology of $\t{B}_n$ with local coefficients;
     \item use of some exact sequences for the cohomology with costant
coefficients.
\end{itemize}

\section{Preliminary results} \label{sec:pre}
In this Section we fix the notation and recall some
preliminary results. We will use classical
 facts (\cite{bourbaki}, \cite{humphreys}) without further reference.
\subsection{Coxeter groups and Artin braid groups}
A \emph{Coxeter graph} is a finite undirected graph, whose edges
are labelled with integers $\geq 3$ or with the symbol $\infty$.

Let  $S$, $E$ be respectively  the vertex and edge set of a
Coxeter graph. For every edge $\{s, t\} \in E$  let $m_{s,t}$  be
its label. If $s, \, t\in S$ ($s \neq t$) are not joined by an
edge, set by convention $m_{s,t}=2$. Let also $m_{s,s}=1$.

Two groups are associated to  a Coxeter graph (as in the
Introduction): the \emph{Coxeter group} $W$ defined by
$$
W=\langle s \in S \, |\, (st)^{m_{s,t}}=1 \,\, \forall s,t \in S
\,\, \textrm{such that}\,\, m_{s,t}\neq \infty\rangle
$$
and the \emph{Artin braid group} $\br{W}$ defined by (see
\cite{Br-Sa}, \cite{brieskorn},  \cite{deligne_imm}):
$$
G=\langle s \in S \, | \,
\begin{underbrace}{stst\ldots}\end{underbrace}_{m_{s,t}-\mathrm{terms}}=\begin{underbrace}{tsts\ldots}\end{underbrace}_{m_{s,t}-\mathrm{terms}}
 \,\, \forall s,t \in S \,\, \textrm{such that}\,\, m_{s,t}\neq
\infty\rangle.
$$
 There is a natural epimorphism $\pi: \br{W}\to W$ and,
by Matsumoto's Lemma \cite{mats}, $\pi$ admits a canonical
set-theoretic section $\psi: W \to \br{W}$.

\subsection{}
In this paper, we are primarily interested in Artin braid groups
associated to Coxeter graphs of type $B_n$, $\t{B}_n$ and $D_n$
(see Table \ref{table:dinkyn}).

\def\objectstyle{\scriptscriptstyle}
\def\mycir {\cir<8pt>{}}

\begin{table}
\begin{center}

\begin{tabular}{ll}
  $B_n$ & \xygraph{
                    !{<0cm, 0cm>; <1.3cm, 0cm>:<0cm, 1.3cm>::}
                    !{(0,0)}*+{1}*\mycir{}="a"
                    !{(1,0)}*+{2}*\mycir{}="b"
                    !{(2,0)}*+{3}*\mycir{}="c"
                    !{(3,0)}*+{4}*\mycir{}="d"
                    !{(3.5,0)}*+{}="va_d"
                    !{(4,0)}*+{}="va_e"
                    !{(4.5,0)}*+{n-2}*\mycir{}="e"
                    !{(5.5,0)}*+{{n-1}}*\mycir{}="f"
                    !{(6.5,0)}*+{{n}}*\mycir{}="g"
                    "a"-"b"
                    "b"-"c"
                    "c"-"d"
                    "d"-@{.}"va_d"
                    "e"-@{.}"va_e"
                    "e"-"f"
                    "f"-"g"^4
         }\\
         {}&{}\\

  $\t{B}_n$ & \xygraph{
                    !{<0cm, 0cm>; <1.3cm, 0cm>:<0cm, 1.3cm>::}
                    !{(0,-0.5)}*+{1}*\mycir{}="a"
                    !{(1,0)}*+{3}*\mycir{}="b"
                    !{(2,0)}*+{4}*\mycir{}="c"
                    !{(2.5,0)}*+{}="va_d"
                    !{(4,0)}*+{}="va_e"
                    !{(4.5,0)}*+{n-1}*\mycir{}="e"
                    !{(5.5,0)}*+{{n}}*\mycir{}="f"
                    !{(6.5,0)}*+{{n+1}}*\mycir{}="g"
                    !{(0,0.5)}*+{{2}}*\mycir{}="z"
                    "a"-"b"
                    "b"-"c"
                    "c"-@{.}"va_d"
                    "e"-@{.}"va_e"
                    "e"-"f"
                    "f"-"g"^4
                    "z"-"b"
         }\\
         {}&{}\\
$D_n$ & \xygraph{
                    !{<0cm, 0cm>; <1.3cm, 0cm>:<0cm, 1.3cm>::}
                    !{(0,-0.5)}*+{1}*\mycir{}="a"
                    !{(1,0)}*+{3}*\mycir{}="b"
                    !{(2,0)}*+{4}*\mycir{}="c"
                    !{(2.5,0)}*+{}="va_d"
                    !{(3,0)}*+{}="va_e"
                    !{(3.5,0)}*+{n-2}*\mycir{}="e"
                    !{(4.5,0)}*+{{n-1}}*\mycir{}="f"
                    !{(5.5,0)}*+{{n}}*\mycir{}="g"
                    !{(0,0.5)}*+{{2}}*\mycir{}="z"
                    "a"-"b"
                    "b"-"c"
                    "c"-@{.}"va_d"
                    "e"-@{.}"va_e"
                    "e"-"f"
                    "f"-"g"
                    "z"-"b"
         }\\

         {}&{}\\
         \end{tabular}
\caption{Coxeter graphs of type $B_n,\; \t{B}_n,\;
D_n$.}\label{table:dinkyn}\end{center}
\end{table}
\def\objectstyle{\textstyle}
The associated Coxeter groups can be described as reflection groups
with respect to an arrangement of hyperplanes (or mirrors). Let $x_1, \ldots, x_n$ be
the standard coordinates in $\R^n$. Consider the
 linear hyperplanes:
\begin{align*}
\mb{H}_k=&\{x_k=0\} &\mb{L}^\pm_{ij}=&\{x_i=\pm x_j\}
\end{align*}
and, for an integer $a \in \Z$, their affine translates:
\begin{align*}
\mb{H}_k(a)=&\{x_k=a\} &\mb{L}^\pm_{ij}(a)=&\{x_i=\pm x_j+a\}
\end{align*}

The Coxeter group $B_n$ is identified with the group of
reflections with respect to the mirrors in the arrangement
$$
\mathcal{A}(B_n):=
                    \{\mb{H}_k\,|\, 1 \leq k \leq n\}
                    \cup
                    \{\mb{L}^\pm_{ij}\,|\, 1 \leq i<j \leq n\}.
$$
As such it is the group of signed permutations of the coordinates
in $\R^n$. Notice that $B_n$ is generated by $n$ basic reflections
$s_1, \ldots, s_n$ having respectively as mirrors
 the $n-1$ hyperplanes $\mb{L}^+_{i,i+1}$ ($1 \leq i \leq n-1$) and the hyperplane $\mb{H}_n$.
 This numbering of the reflections is consistent with the numbering
 of the vertices of the Coxeter graph for $B_n$ shown in Table \ref{table:dinkyn}.

The affine Coxeter group $\t{B}_n$ is the semidirect product of the
Coxeter group $B_n$ and the coroot lattice, consisting of integer
vectors whose coordinates add up to an even number. The arrangement
of mirrors is then the affine hyperplane arrangement:
\begin{equation}\label{e:refl:tb}
\mathcal{A}(\t{B}_n):=
                    \{\mb{H}_k(a)\, |\, 1\leq k \leq n, \, a \in \Z\}
                    \cup
                    \{\mb{L}^\pm_{ij}(a)\, |\, 1 \leq i <j \leq n, \, a \in \Z\}.
\end{equation}
It is generated by the basic reflections for $B_n$ plus an extra
affine reflection $\t{s}$ having $\mb{L}^-_{12}(1)$ as mirror. The latter
commutes with all the basic reflections of $B_n$ but $s_2$, for which
$(\t{s}{s_2})^3=1$. This accounts for the Coxeter graph of type
$\t{B}_n$ in the table, where, however, we chose by our convenience
a somewhat unusual vertex numbering.

Finally the group $D_n$ has reflection arrangement:
$$
\mathcal{A}(D_n):=\{\mb{L}^\pm_{ij}\, | \, 1 \leq i < j \leq n\}
$$
and it can be regarded as the group of signed permutations of the
coordinates which involve an even number of sign changes. In
particular $D_n$ is a subgroup of index $2$ in $B_n$. The group is
generated by $n$ basic reflections w.r.t. the hyperplanes
$\mb{L}^-_{12}$ and $\mb{L}^+_{i,i+1}$ ($1 \leq i \leq n-1$).

\subsection{Generalized Poincar\'e series}\label{sec:poincare}
For future use in cohomology computations, we will need some
analog of ordinary Poincar\'e series for Coxeter groups. Consider
a domain $R$ and let $R^*$ be the group of unit of $R$. Given an
abelian representation
$$
\eta:\, \br{W} \to R^*
$$
of the Artin group $\br{W}$ and a finite subset ${U}\sst W$, we may
consider the $\eta$-Poincar\'e series:
$$
U(\eta)=\sum_{w \in U}(-1)^{\ell(w)} \eta(\psi w)\in R
$$
where $\ell$ is the length in the Coxeter group and $\psi: W \to
\br{W}$ is the canonical section. In particular, when $W$ is
finite, we say that $W(\eta)$ is the $\eta$-Poincar\'e series of
the group. Notice that for $R=\Q[q^\pmu]$ we may consider the
representation  $\eta_q$ that sends the standard generators of
$\br{W}$ into $(-q)$-multiplication; in this situation we recover
the ordinary Poincar\'e series:
$$
W(\eta_q)=W(q)
$$
Further, for the Artin group of type $W=B_n, \, \t{B}_n$ we are
interested in the representation
$$\eta_{q,t}: \, \br{W} \to \Q[q^\pmu, t^\pmu]$$
defined sending the last standard generator (the one laying in the
tree leave labelled with $4$) to $(-t)$-multiplication and the
remaining ones to $(-q)$-multiplication. The associated Poincar\'e
series $B_n(q,t):=B_n(\eta_{q,t})$ will be called the
\emph{$(q,t)$-weighted Poincar\'e series} for $B_n$.

In order to recall closed formulas for Poincar\'e series, we first
 fix some notations that will be adopted throughout the paper.
We define the $q$-analog of a positive integer $m$ to be the
polynomial
$$[m ]_q := 1 + q + \cdots q^{m-1} = \frac{q^m -1}{q-1}$$
It is easy to see that $[m ] = \prod_{i \mid m} \ph_m(q)$.
Moreover we define the $q$-factorial and double factorial
inductively as:
\begin{flalign*}
[m]_q!&:=[m]_q \cdot [m-1]_q!\\
[m]_q!!&:=[m]_q \cdot [m-2]_q!!
\end{flalign*}
where is understood that $[1]!=[1]!!=[1]$ and $[2]!!=[2]$. A
$q$-analog of the binomial $\binom{m}{i}$ is given by the
polynomial
$$\qbin{m}{i}_{q}: = \frac{[m]_q!}{[i]_q! [m-i]_q!}$$
We can also define the $(q,t)$-analog of an even number
$$[2m]_{q,t} := [m]_q (1+tq^{m-1})$$
and of the double factorial
$$
[2m]_{q,t}!! := \prod_{i=1}^m [2i]_{q,t}\ =\ [ m ]_q!
\prod_{i=0}^{m-1} (1+tq^i)
$$
Notice that specializing $t$ to $q$, we recover the $q$-analogue of
an even number and of its double factorial. Finally, we define the
polynomial
\begin{equation}\label{e:binqt}
\qbin{m}{i}_{q,t}': = \frac{[2m]_{q,t}!!}{[2i]_{q,t}!! [m-i]_q!} \
= \qbin{m}{i}_{q}\prod_{j=i}^{m-1}(1+tq^j)
\end{equation}
With this notation the ordinary Poincar\'e series for $D_n$ and
$B_n$ may be written as
\begin{align}
    D_n(q)&:=\sum_{w \in D_n}q^{\ell(w)}= [2(n-1)]_q!! \cdot [n]_q \label{e:d:q:poincare}\\
    B_n(q)&:=\sum_{w \in B_n}q^{\ell(w)}=[2n]_q!! \label{e:b:q:poincare}
\end{align}
while the \emph{$(q,t)$-weighted Poincar\'e series} for $B_n$ is
given by (see e.g. \cite{reiner}):
\begin{equation} \label{p:qtpoincare}
B_n(q,t)=
[2n]_{q,t}!!
\end{equation}

\section{The $K(\pi, 1)$ problem for the affine Artin group of type $\t{B}_n$} \label{sec:kpi_1}

Using the explicit description of the reflection mirrors in Equation (\ref{e:refl:tb}), the complement of the complexified affine
reflection arrangement of type $\t{B}_n$ is given by:
$$
\mb{Y}:=\mb{Y}(\t{B}_n)=\{ x \in \C^n \, | \, x_i \pm x_j \notin \Z \textrm{
for all $i\neq j$}, \, x_k \notin \Z \textrm{ for all $k$} \}
$$
On $\mb{Y}$ we have, by standard facts,  a free action by translations
of the coweight lattice $\Lambda$, identified with the standard
lattice $\Z^n\subset \C^n$.
\bigskip

 \noindent {\bf{Proof of Theorem \ref{thm:main}}}\
We first  explicitly describe the covering $\mb{Y} \to \mb{Y}/\Lambda$
applying the exponential map $y=\exp(2 \pi \i x)$ componentwise to
$\mb{Y}$:
\[
\xymatrix@R=12pt{
  \mb{Y}  \ar[r]^-\pi &  \mb{Y}/\Lambda\cong\{y \in \C^n \, | \, y_i \neq y_j^{\pmu},\,\, y_k \neq 0,1\}\\
 (x_1, \ldots, x_n) \ar@{|->}[r] & \big(\exp{(2 \pi \i x_1)}, \ldots, \exp{(2 \pi \i x_n)}\big)}
\]
Notice now that the function
$$
\C\setminus\{0,1\} \ni y \mapsto g(y)= \frac{1+y}{1-y} \in
\C\setminus\{\pm 1\}
$$
satisfies $g(y^{-1})=-g(y)$. Further $g$ is invertible, its
inverse being given by $z \mapsto \frac{z-1}{z+1}$. Therefore
applying $g$ componentwise to $\mb{Y}/\Lambda$, we have:
$$
\mb{Y}/\Lambda \cong \{z \in \C^n\, | \, z_i \neq \pm z_j, \,\, z_k
\neq \pmu\}
$$

Consider now the arrangement $\mathcal{A}$ in $\R^{n+1}$
consisting of the hyperplanes $\mb{L}^\pm_{ij}$ for $1\leq i<j \leq
n+1$ and $\mb{H}_{1}$ and let $\mb{Y}(\mathcal{A})$ be the complement of its
complexification.

We have an homeomorphism $$\eta: \C^* \times \mb{Y}/\Lambda \to
\mb{Y}(\mathcal{A})$$ defined by $$\eta\big(\lambda, (z_1, \ldots, z_n)
\big)=(\lambda, \lambda z_1, \ldots, \lambda z_n)$$

 To show that $\mb{Y}/\Lambda$ is a $K(\pi,1)$, it is then sufficient to show that
$\mb{Y}(\mathcal{A})$ is a $K(\pi,1)$.
 We will show in Lemma \ref{lemma:simplicial} below that $\mathcal{A}$
 is simplicial, and therefore the result follows from Deligne's
 Theorem \ref{thm:deligne}.  \qed
\vs
 \noindent {\bf{Remark}}\ By the same exponential argument one may recover the results of \cite{okonek} for the affine Artin group of type $\t{A}_n, \t{C}_n$
 (for further applications we refer to \cite{allcock-braid}).

\begin{lemma}\label{lemma:simplicial}
Let $\mathcal{A}$ be the real arrangement in $\R^{n+1}$ consisting
of the hyperplanes $\mb{L}^\pm_{ij} $ for $1\leq i<j \leq n+1$ and
$\mb{H}_1$. Then $\mathcal{A}$  is simplicial.
\end{lemma}
\begin{dm}
Notice that $\mathcal{A}$ is the union of the reflection arrangement
$\mathcal{A}(D_{n+1})$ of type $D_{n+1}$ and the hyperplane
$\mb{H}_{1}=\{x_{1}=0\}$. Hence we study how the chambers of
$\mathcal{A}(D_{n+1})$ are cut by the hyperplane $\mb{H}_1$. Since the
Coxeter group $D_{n+1}$ acts transitively on the collection of
chambers, it is enough to consider how the fundamental chamber
$\mb{C}_0$ of $\mathcal{A}(D_{n+1})$ is cut by the $D_{n+1}$-translates
of the hyperplane $\mb{H}_1$, i.e.~by the coordinate hyperplanes $\mb{H}_k$
for $k=1,2,\ldots, n+1$.\\
We may choose
$$
\mb{C}_0=\{-x_2<x_1<x_2< \ldots <x_n<x_{n+1}\}
$$
as fundamental chamber. Of course, this is a simplicial cone. Notice
that the coordinate of a point in $\mb{C}_0$ are all positive except
(possibly) the first. Thus it is clear that for $k
\geq 2 $ the hyperplanes $\mb{H}_k$ do not cut $\mb{C}_0$.\\
A quick check shows instead that $\mb{H}_1$ cuts $\mb{C}_0$ into two
simplicial cones $\mb{C}_1$, $\mb{C}_2$ given precisely by:
\begin{align*}
\mb{C}_1&=\{0<x_1<x_2< \ldots <x_n<x_{n+1}\}\\
\mb{C}_2&=\{0<-x_1<x_2< \ldots <x_n<x_{n+1}\}
\end{align*}
\qed
\end{dm}


\section{Cohomology} \label{sec:cohom}
In this Section we will compute the cohomology groups
$$H^*(G_{\widetilde{B}_n}, \Q[[q^\pmu, t^\pmu]]_{q,t})$$ where $\Q[[q^\pmu, t^\pmu]]_{q,t}$ is the local system over the module of
Laurent series $\Q [[q^\pmu,
t^\pmu]]$ and the action is $(-q)-$multiplication for the standard
generators associated to the first $n$ nodes of the Dynkin
diagram, while is $(-t)-$multiplication for the generator
associated to the last node.

\subsection{Algebraic complexes for Artin groups}

As a main tool for  cohomological computations we use the
algebraic complex described in \cite{salvettiArtin} (see the Introduction);
the algebraic generalization of this complex by De Concini-Salvetti
\cite{salvettiDeconciniNotaArtin} provides an effective way to
determine the cohomology of the orbit space $X(W)$ with values in
an arbitrary $\br{W}$-module. When $X(W)$ is a $K(\pi,1)$ space,
of course, we get the cohomology of the group $\br{W}$.

For sake of simplicity, we restrict ourself to the abelian
representations considered in Section \ref{sec:poincare}. Let
$(W,S)$ be a Coxeter system. Given a a representation
$\eta:\br{W}\to R^*$, let $M_\eta$ be the induced structure of
$\br{W}$-module on the $R$-module $M$. We may describe a cochain complex $C^*(W)$
for the cohomology $H^*(X(W); M_\eta)$ as follows. The cochains in
dimension $k$ consist in the free $R$-module indexed by the finite
parabolic subgroup of $W$:
\begin{equation}\label{eq:chain:compl:braids}
C^k(W):= \bigoplus_{
    \begin{array}{c}
        \scriptstyle {\Gamma}: \card{\Gamma}=k \\
        \scriptstyle \card{W_{\Gamma}}< \infty
    \end{array}}
    M.e_\Gamma
\end{equation}
and the coboundary map are completely described by the formula:
\begin{equation}\label{eq:chain:bordo:braids}
\di (e_\Gamma)=   \sum_{
                             \begin{array}{c}
                                \scriptstyle \Gamma'\tss \Gamma\\
                                \scriptstyle \card{\Gamma'}=\card{\Gamma}+1\\
                                \scriptstyle \card{W_{\Gamma'}}< \infty
                            \end{array}
                            }
                            (-1)^{\alpha({\Gamma}, \Gamma')}
                            \frac{W_{\Gamma'}(\eta)}{W_{\Gamma}(\eta)}
                             e_{\Gamma'}
\end{equation}
where $W_\Gamma(\eta)$ is the $\eta$-Poincar\'e series of the
parabolic subgroup $W_\Gamma$ and ${\alpha(\Gamma, \Gamma')}$ is
an incidence index depending on a fixed linear order of $S$. For
$\Gamma'\setminus \Gamma=\{s'\}$ it is defined as
$$
    {\alpha({\Gamma}, \Gamma')}:= \card{\{s \in \Gamma \,:\ \, s< s'\} }
$$

 We identify (consistently with  Table \ref{table:dinkyn}) the
generating reflections set $S$ for $\t{B}_n$ with the set $\{1,2,
\ldots, n+1\}$.
It is useful to represent a subset $\Gamma \sst
S$ with its characteristic function. For example the subset
$\{1,3,5,6\}$ for $\t{B}_6$ may be represented as the binary
string:
$$
 \abin{0}{1}10110
$$
To determine the cohomology of $\br{\t{B}_n}$, it will be
necessary to give a close look to the cohomology of $\br{D_n}$. It
is convenient to number the vertex of $D_n$ as in table
\ref{table:dinkyn} and to regard  parabolic subgroups as
binary strings as before. 

\subsection{}\label{par:shift}

Let $R$ be the ring of Laurent polynomials $\Q[q^\pmu, t^\pmu]$ and $M$
be the $R$-module of Laurent series $\Q[[q^\pmu, t^\pmu]]$ and let
$R_{q,t}$, $M_{q,t}$ be the corresponding local systems, with
action $\eta_{q,t} $. Our main interest is to compute the
cohomology with trivial rational coefficient of the group
$$
Z_{\widetilde{B}_n}=\ker {(G_{\widetilde{B}_n} \to \Z^2)}
$$
that is the commutator subgroup of $G_{\widetilde{B}_n}$.
By Shapiro Lemma (see \cite{brown}) we have the following equivalence:
$$
H^*(Z_{\widetilde{B}_n}, \Q) \simeq H^*(G_{\widetilde{B}_n}, M_{q,t})
$$
and the second term of the equality is computed by the Salvetti
complex $C^*(\widetilde{B}_n)$ over the module $M_{q,t}$. Notice
that the finite parabolic subgroups of $W_{\widetilde{B}_n}$ are
in $1-1$
correspondence with the proper subsets of the set of simple roots
$S$.
%

We can define an \emph{augmented} Salvetti complex $\widehat{C}^*(\widetilde{B}_n)$ as follows:
$$
\widehat{C}^*(\widetilde{B}_n) = C^*(\widetilde{B}_n) \oplus {(M_{q,t})}.e_{S}.
$$
We need to define the boundary map for the $n$-dimensional
generators. Let we first define a quasi-Poincar\'e polynomial for
$G_{\widetilde{B}_n}$. We set
$$
\widehat{W}_{S}(q,t) = \widehat{W}_{\widetilde{B}_n}(q,t) =
[2(n-1)]!!\ [n]\ \prod_{i=0}^{n-1} (1+tq^i) .
$$
It is easy to verify that $\widehat{W}_{\widetilde{B}_n}(q,t)$ is the least common multiple of all $W_\Gamma(q,t)$, for $\Gamma \sst S$ with $|\Gamma| = n$. This allows us to define the boundary map for the generators $e_\Gamma$, with $|\Gamma| = n$:
$$
d(e_\Gamma) = (-1)^{\alpha(\Gamma, S)} \frac{\widehat{W}_{\widetilde{B}_n}(q,t)}{W_\Gamma(q,t)}e_{S}
$$
and it is straightforward to verify that $\widehat{C}^*(\widetilde{B}_n)$ is still a chain complex.
Moreover we have the following relations between the cohomologies of $C^*(\widetilde{B}_n)$ and $\widehat{C}^*(\widetilde{B}_n)$:
$$
H^i(C^*(\widetilde{B}_n)) = H^i(\widehat{C}^*(\widetilde{B}_n))
$$
for $i \neq n, n+1$ and we have the short exact sequence
$$
0 \to H^n(\widehat{C}^*(\widetilde{B}_n), M_{q,t}) \to H^n(C^*(\widetilde{B}_n), M_{q,t}) \to M_{q,t} \to 0.
$$
Finally one can prove that the complex $\widehat{C}^*(\widetilde{B}_n)$ with coefficients in
the local system $R_{q,t}$ is \emph{well
filtered} (as defined in \cite{C}) with respect to the variable
$t$ and so it gives
%
%
%
%
the same cohomology, modulo an index shifting, of the complex with
coefficients over the module $\Q[t^\pmu][[q^\pmu]]$. Another index
%
%
%
shifting can be proved with a slight improvement of the results in
\cite{C}, allowing to pass to the module $M$. Hence we have the
following
\begin{prop}
$$
H^i(Z_{\widetilde{B}_n}, \Q) \simeq
H^i(\widehat{C}^*(\widetilde{B}_n), M_{q,t})\simeq
H^{i+2}(\widehat{C}^*(\widetilde{B}_n), R_{q,t}) \simeq
H^{i+2}(G_{\widetilde{B}_n}, R_{q,t})
$$
for $i \neq n, n+1$ and
$$
H^{n}(Z_{\widetilde{B}_n}, \Q) \simeq H^{n}(G_{\widetilde{B}_n}, M_{q,t}) \simeq M
$$
$$
H^{n+1}(Z_{\widetilde{B}_n}, \Q) \simeq H^{n+1}(G_{\widetilde{B}_n}, M_{q,t})\simeq 0.
$$
\qed
\end{prop}

From now on we deal only with the complex
$\widehat{C}^*(\widetilde{B}_n)$ with coefficients in the local
system $R_{q,t}$.

\subsection{\label{s:splitting}}
For Coxeter groups of type $W=D_n$, $\t{B}_n$  the Salvetti's
complex $C^*W$ exhibits an involution $\sigma$ defined by:
\begin{align*}
    \abin{0}{0}A&\stackrel{\sigma}{\longrightarrow} \abin{0}{0}A &
    \abin{1}{1}A&\stackrel{\sigma}{\longrightarrow}-\abin{1}{1}A\\
    \abin{0}{1}A&\stackrel{\sigma}{\longrightarrow}\abin{1}{0}A &
    \abin{1}{0}A&\stackrel{\sigma}{\longrightarrow}\phantom{-}\abin{0}{1}A.
\end{align*}
Let $\inv^* W $ be the module of $\sigma$-invariants and
$\anti^*W$ the module of $\sigma$-anti-invariants. We may then
split the complex into:
\begin{align*}
C^* W &= \inv^* W \oplus \anti^* W.
\end{align*}
In particular the computation of the cohomolgy of $C^*W$ may be
performed analyzing separately the two subcomplexes.

\subsection{Cohomology of $\anti^*D_n$} \label{ss:anti-inv}

The cohomology of the anti-invariant subcomplex for $D_n$ was
completely determined in \cite{arithmetic_artin}. However we will
need for our purposes generators for the cohomology groups which
are not easily deduced from the argument in the original paper. So
we briefly recall this result.

Let $G^1_n$ be the subcomplex of $C(D_n)$ generated by the strings
of type $\abin{0}{1}A$ and $\abin{1}{1}A$. It is easy to see that
$G^1_n$ is isomorphic (as a complex) to $\anti(D_n)$.

Define the set
$$
S_n = \{ h\in \N \mbox{ s. t. }2h | n \mbox{ or } h|n-1 \mbox{ and
} 2h \nmid (n-1)\}
$$
Note that $h$ appears in $S_n$ if and only if $n=2\lambda h$
(i.e.~$n$ is an even multiple of $h$) or $n=(2\lambda+1)h +1$ ($n$
is an odd multiple of $h$ incremented by $1$).

\begin{prop}[\cite{arithmetic_artin}] \label{p:dpssG}
The top-cohomology of $G^1_n$ is:
$$
H^n G^1_n = \bigoplus_{h \in S_n}\{2h\}
$$
whereas for $s>0$ one has:
\begin{align*}
H^{n-2s}G^1_n &=
                 \bigoplus_{
                        \begin{array}{c}
                        \scriptstyle
                             h \in S_n \\
                            1 < h < \frac{n}{2s}
                        \end{array} }
                 \{ 2h \}\\
H^{n-2s+1}G^1_n &=
                \bigoplus_{
                        \begin{array}{c}
                        \scriptstyle
                            h \in S_n \\
                            1 < h \leq \frac{n}{2s}
                        \end{array} }
                \{ 2h \}.
\end{align*} \qed
\end{prop}
\vs
We need a description of the generators for these modules.\\
First we define the following basic binary strings:
\begin{align*}
o_{\mu}[h]&=\left\{ \begin{array}{ll}
                            \abin{0}{1}1^{h-1} & \textrm{    for
                            $\mu=0$}\\ {}&\\
                            \abin{1}{1}1^{2\mu h -2}01^h & \textrm{    for $\mu \geq
                            1$}
                    \end{array}
            \right. \\
e_{\mu}[h]&=\abin{1}{1}1^{(2\mu -1)h -1}01^{h-2} \mbox{   for }\mu
\geq 1
\end{align*}
\begin{align*}
 s_h&= 01^{h-2} &
 l_h&= 01^h.
\end{align*}
A set of candidate cohomology generators is given by the following
cocycles:
\begin{align*}
    o_{\mu, 2i}[h]&= \frac{1}{\ph_{2h}}d(o_{\mu}[h](s_h l_h)^i)\\
    o_{\mu, 2i+1}[h]&= \frac{1}{\ph_{2h}}d(o_{\mu}[h](s_h
        l_h)^is_h)\\
    e_{\mu, 2i}[h]&= \frac{1}{\ph_{2h}}d(e_{\mu}[h](l_h s_h)^i)\\
    e_{\mu, 2i+1}[h]& = \frac{1}{\ph_{2h}}d(e_{\mu}[h](l_h s_h)^il_h).
\end{align*}
Indeed these cocycles account for all the generators:
\begin{prop}\label{p:anti:gen}
\begin{enumerate}
    \item Let $n=2\lambda h$. Then for $0\leq s < \lambda$ the summand
        of $H^{n-2s}(G^1_n)$ isomorphic to $\{2h\}$ is generated by
        $e_{\lambda-s,2s}[h]$. Similarly for $0\leq s < \lambda$ the
        summand of $H^{n-2s-1}(G^1_n)$ is generated by
        $o_{\lambda-s-1,2s+1}[h]$.
    \item Let $n=(2\lambda +1)h+1$. Then for $0\leq s \leq \lambda$ the summand of $H^{n-2s}(G^1_n)$ isomorphic to
        $\{2h\}$ is generated by $o_{\lambda-s,2s}[h]$. For $0\leq s <
        \lambda$ the summand of $H^{n-2s-1}(G^1_n)$ is generated by
        $e_{\lambda-s,2s+1}[h]$.
\end{enumerate}
\end{prop}
Proposition \ref{p:anti:gen} is best proven by induction on $n$,
recovering in particular the quoted result from
\cite{arithmetic_artin}.
\begin{dm}
We filter the complex $G^1_n$ from the right and use the
associated spectral sequence. Let:
$$
F_k G^1_n= \langle A 1^k \rangle
$$
be the subcomplex generated by binary strings ending with at least
$k$ ones. We have a filtration
$$
G^1_n = F_0 G^1_n \tss F_1 G^1_n \tss \ldots \tss F_{n-2} G^1_n
\tss F_{n-1} G^1_{n-1}\tss 0
$$
 in which the subsequent quotients for $k=1,2, \ldots, n-3$
$$
\frac{F_k G^1_n}{  F_{k+1} G^1_n} = \langle A 0 1^k\rangle \cong
G^1_{n-k-1}[k]
$$
are isomorphic to the complex for $G^1_{n-k-1}$ shifted in degree
by $k$, while
\begin{align*}
    \frac{F_{n-2} G^1_n}{  F_{n-1} G^1_n} = & \left\langle
        \abin{0}{1} 1^{n-2} \right\rangle \cong R[n-1] &
    {F_{n-1} G^1_n} = & \left\langle
        \abin{1}{1} 1^{n-2} \right\rangle \cong R[n]. &
\end{align*}
Therefore the columns of the $E_1$ term of the spectral sequence
are either the module $R$ or  are given by the cohomology of
$G^1_{n'}$ with $n'<n$. Reasoning by induction, we may thus
suppose that their cohomology has the generators prescribed by the
proposition. Since there can be no non-zero maps between the
module $\{2h\}$, $\{2h'\}$ for $h \neq h'$, we may separately
detect the $\ph_{2h}$-torsion in the cohomology.\\
Fix an integer $h>1$. Then the relevant modules for the
$\ph_{2h}$-torsion in the $E_1$ term are suggested in Table
\ref{table:e1:anti}. We will call a column \emph{even} if it is
relative to $G^1_{2\mu h}$ and \emph{odd} if it is relative to
$G^1_{(2\mu+1 )h+1}$ for some $\mu$.
\begin{table}[ht!]
\begin{center}
\def\objectstyle{\scriptstyle}
\def\myob {\phantom{o_{a,c}}}
$$
\xygraph{  
                    !{<0cm, 0cm>; <0cm, .4cm>:<.65cm, 0cm>::}
                    !{(16,0)}*+{\myob}*\frm{-,}="5-0"    !{(16,0)}*+{o_{2,0}}="5-0b"
                    !{(15,0)}*+{\myob}*\frm{-,}="5-1"    !{(15,0)}*+{e_{2,1}}="5-1b"
                    !{(14,0)}*+{\myob}*\frm{-,}="5-2"    !{(14,0)}*+{o_{1,2}}="5-2b"
                    !{(13,0)}*+{\myob}*\frm{-,}="5-3"    !{(13,0)}*+{e_{1,3}}="5-3b"
                    !{(12,0)}*+{\myob}*\frm{-,}="5-4"    !{(12,0)}*+{o_{0,4}}="5-4b"
                    !{(12,4)}*+{\myob}*\frm{-,}="4-0"    !{(12,4)}*+{e_{2,0}}="4-0b"
                    !{(11,4)}*+{\myob}*\frm{-,}="4-1"    !{(11,4)}*+{o_{1,1}}="4-1b"
                    !{(10,4)}*+{\myob}*\frm{-,}="4-2"    !{(10,4)}*+{e_{1,2}}="4-2b"
                    !{(9,4)}*+{\myob}*\frm{-,}="4-3"     !{(9,4)}*+{o_{0,3}}="4-3b"
                    !{(10,6)}*+{\myob}*\frm{-,}="3-0"    !{(10,6)}*+{o_{1,0}}="3-0b"
                    !{(9,6)}*+{\myob}*\frm{-,}="3-1"     !{(9,6)}*+{e_{1,1}}="3-1b"
                    !{(8,6)}*+{\myob}*\frm{-,}="3-2"     !{(8,6)}*+{o_{0,2}}="3-2b"
                    !{(6,10)}*+{\myob}*\frm{-,}="2-0"    !{(6,10)}*+{e_{1,0}}="2-0b"
                    !{(5,10)}*+{\myob}*\frm{-,}="2-1"    !{(5,10)}*+{o_{0,1}}="2-1b"
                    !{(4,12)}*+{\myob}*\frm{-,}="1-0"    !{(4,12)}*+{o_{0,0}}="1-0b"
                    !{(2,14)}*+{R}="libb"
                    !{(2,15)}*+{R}="topb"
                    !{(16,-1)}*+{}="y"
                    !{(0,16)}*+{}="x"
                    !{(0,-2)}*+{}="o1"
                    !{(-1,-1)}*+{}="o2"
                     !{(-1,0)}*+{G^1_{5h+1}}="C5"
                    !{(-1,4)}*+{G^1_{4h}}="C4"
                    !{(-1,6)}*+{G^1_{3h+1}}="C3"
                    !{(-1,10)}*+{G^1_{2h}}="C2"
                    !{(-1,12)}*+{G^1_{h+1}}="C1"
                    "5-1":@/^/^{d_{h+1}}"4-0"
                    "5-2":@/^/"4-1"
                    "5-3":@/^/"4-2"
                    "5-4":@/^/"4-3"
                    "4-1":@/_/"3-0"
                    "4-2":@/_/"3-1"
                    "4-3":@/_/_{d_{h-1}}"3-2"
                    "3-1":@/^/^{d_{h+1}}"2-0"
                    "3-2":@/^/"2-1"
                    "2-1":@/_/_{d_{h-1}}"1-0"
                    "libb":_{d_1}"topb"
                    "o1":"x"  
                    "o2":"y"  
}
$$
\caption{Spectral sequence for $G^1_n$}\label{table:e1:anti}
\end{center}
\end{table}
The differential $d_1$ is zero everywhere but $d_1:\,
E_1^{(n-2,1)}\to E_1^{(n-1,1)} $ where it is given by
multiplication by $[2(n-1)]!!/[n-1]!$. Thus the $E_2$ term differs
from the $E_1$ only in positions $(n-2,1)$ and $(n-1,1)$, where:
\begin{align*}
    E_2^{(n-2,1)}&=0 &
    E_2^{(n-1,1)}&=\frac{R}{[2(n-1)]!!/[n-1]!}
\end{align*}
Then all other differentials are zero up to $d_{h-2}$.\\
 It is now useful to distinguish among 4 cases according to the
remainder of $n \, \mathrm{mod} ( 2h )$:
\begin{itemize}
    \item[a)] $n=2 \lambda h +c$ for $1 \leq c \leq h$
    \item[b)] $n=(2 \lambda+1) h + 1$
    \item[c)] $n=(2 \lambda +1)h +1+ c $ for $1 \leq c \leq h-2$
    \item[d)] $n=2 \lambda h$-
\end{itemize}
In case a), note the first column relevant for $\ph_{2h}$-torsion
is even (see also Table \ref{table:eh-1:anti:case:a}).

\begin{table}[ht]
\begin{center}
\def\objectstyle{\scriptstyle}
\UsePatternFile{macpat.xyp}
    \AliasPattern{bricks:c}{mac14}{xymacpat} 
    \AliasPattern{bricks}{mac06}{xymacpat}
\def\myob {\phantom{o_{a,c}}}
$$
\xygraph{  
                    !{<0cm, 0cm>; <0cm, .25cm>:<.6cm, 0cm>::}
%
                    !{(16,0)}*+{}*[bricks]\frm{**}="8-0"
                    !{(15,0)}*+{}*[bricks]\frm{**}="8-1"
                    !{(14,0)}*+{}*[bricks]\frm{**}="8-2"
                    !{(13,0)}*+{}*[bricks]\frm{**}="8-3"
                    !{(12,0)}*+{}*[bricks]\frm{**}="8-4"
                    !{(11,0)}*+{}*[bricks]\frm{**}="8-5"
                    !{(10,0)}*+{}*[bricks]\frm{**}="8-6"
                    !{(9,0)}*+{}*[bricks]\frm{**}="8-7"
                    !{(14,2)}*+{}*[bricks]\frm{**}="7-0"
                    !{(13,2)}*+{}*[bricks]\frm{**}="7-1"
                    !{(12,2)}*+{}*[bricks]\frm{**}="7-2"
                    !{(11,2)}*+{}*[bricks]\frm{**}="7-3"
                    !{(10,2)}*+{}*[bricks]\frm{**}="7-4"
                    !{(9,2)}*+{}*[bricks]\frm{**}="7-5"
                    !{(8,2)}*+{}*[bricks]\frm{**}="7-6"
                    !{(12,4)}*+{\dots\dots}="dots2"
                    !{(11,4)}*+{\dots\dots}="dots3"
                    !{(10,4)}*+{\dots\dots}="dots4"
                    !{(9,4)}*+{\dots\dots}="dots5"
                    !{(8,4)}*+{\dots\dots}="dots6"
                    !{(7,4)}*+{\dots\dots}="dots7"
                    !{(10,6)}*+{}*[bricks]\frm{**}="4-0"
                    !{(9,6)}*+{}*[bricks]\frm{**}="4-1"
                    !{(8,6)}*+{}*[bricks]\frm{**}="4-2"
                    !{(7,6)}*+{}*[bricks]\frm{**}="4-3"
                    !{(8,8)}*+{}*[bricks]\frm{**}="3-0"
                    !{(7,8)}*+{}*[bricks]\frm{**}="3-1"
                    !{(6,8)}*+{}*[bricks]\frm{**}="3-2"
                    !{(4,12)}*+{}*[bricks]\frm{**}="2-0"
                    !{(3,12)}*+{}*[bricks]\frm{**}="2-1"
                    !{(2,14)}*+{}*[bricks]\frm{**}="1-0"
                    !{(0,17)}*+{\frac{R}{(\ph_{2h})^\lambda}}="topb"
                    !{(16,-2)}*+{}="y"
                    !{(-2,18)}*+{}="x"
                    !{(-2,-2.5)}*+{}="o1"
                    !{(-3,-2)}*+{}="o2"
                     !{(-3,0)}*+{G^1_{2 \lambda h}}="Cpari"
                    !{(-3,2)}*+{G^1_{(2 \lambda -1)h +1}}="Cdisp"
                    !{(-3,6)}*+{G^1_{4h}}="C4"
                    !{(-3,8)}*+{G^1_{3h+1}}="C3"
                    !{(-3,12)}*+{G^1_{2h}}="C2"
                    !{(-3,14)}*+{G^1_{h+1}}="C2"
                    "o1":"x"  
                    "o2":"y"  
                     "8-1":"7-0"
                     "8-2":"7-1"
                     "8-3":"7-2"
                     "8-4":"7-3"
                     "8-5":"7-4"
                     "8-6":"7-5"
                     "8-7":_{d_{h-1}}"7-6"
                     "4-1":"3-0"
                     "4-2":"3-1"
                     "4-3":_{d_{h-1}}"3-2"
                     "2-1":_{d_{h-1}}"1-0"
}
$$
\caption{$E_{h-1}$-term of the spectral sequence for $G^1_n$ in
case a)}\label{table:eh-1:anti:case:a}
\end{center}
\end{table}

The differential $d_{h-1}$ maps the modules of positive
codimension of an even column $G^1_{2\mu h}$ ($1 \leq \mu \leq
\lambda$)
  to those in the odd column $G^1_{(2\mu -1)h+1}$. Using
the suitable generators of type $e_{\cdot, \cdot}[h], o_{\cdot,
\cdot}[h]$, the map $d_{h-1}$ may be identified with the
multiplication by
\begin{equation}\label{eq:diff:h-1}
\qbin{n-(2\mu-1)h-1}{h-1}=\qbin{2(\lambda-\mu)+c+h-1}{h-1}
\end{equation}
Since this polynomial is non-divisible by $\ph_{2h}$, the
restriction of $d_{h-1}$ to positive codimension elements in even
columns  is injective. It follows that in the $E_h$-term the only
survivors are in positions $(c+ 2( \lambda- \mu)h-1, 2\mu h)$,
generated by $e_{\mu, 0}[h]$ and
$$
E_h^{(n-1,1)}\cong E_2^{(n-1,1)}=\frac{R}{[2(n-1)]!!/[n-1]!}.
$$
Note that  in $E_h^{(n-1,1)}$ the only  torsion of type
$\ph_{2h}^l$ is given by the summand:
$$
\frac{R}{(\ph_{2h})^\lambda}
$$
The setup is summarized in Table \ref{table:ehigh:anti:case:a}. In
the Table the survivors are in dark grey boxes while annihilated
terms are in light grey.

Further, using the generators and up to an invertible, we may
identify the differential $d_{2 \mu h}: E_{2\mu h}^{(c+ 2(
\lambda- \mu)h-1, 2\mu h)} \to E_{2\mu h}^{n-1,1} $ with the
multiplication by $\ph_{2h}^{\lambda-\mu}$ ($1 \leq \mu \leq
\lambda$). Thus, for example, in the $E_{2h+1}$ term the module in
position $(c+ 2( \lambda- 1)h-1, 2 h)$ vanishes and the
$\ph_{2h}$-torsion in $E_{2h+1}^{(n-1,1)}$ is reduced to
${R}/{(\ph_{2h})^{\lambda-1}}$. Continuing in this way, all
$\ph_{2h}$-torsion vanishes. In summary there is no
$\ph_{2h}$-torsion in the cohomology of $G^1_n$; this ends case
a).
\begin{table}[ht]
\begin{center}
\def\objectstyle{\scriptstyle}
\UsePatternFile{macpat.xyp}
    \AliasPattern{bricks:c}{mac14}{xymacpat} 
    \AliasPattern{bricks}{mac06}{xymacpat}
\def\myob {\phantom{o_{a,c}}}
$$
\xygraph{  
                    !{<0cm, 0cm>; <0cm, .25cm>:<.6cm, 0cm>::}
%
                    !{(16,0)}*+{}*[bricks]\frm{**}="8-0"
                    !{(15,0)}*+{}*[bricks:c]\frm{**}="8-1"
                    !{(14,0)}*+{}*[bricks:c]\frm{**}="8-2"
                    !{(13,0)}*+{}*[bricks:c]\frm{**}="8-3"
                    !{(12,0)}*+{}*[bricks:c]\frm{**}="8-4"
                    !{(11,0)}*+{}*[bricks:c]\frm{**}="8-5"
                    !{(10,0)}*+{}*[bricks:c]\frm{**}="8-6"
                    !{(9,0)}*+{}*[bricks:c]\frm{**}="8-7"
                    !{(14,2)}*+{}*[bricks:c]\frm{**}="7-0"
                    !{(13,2)}*+{}*[bricks:c]\frm{**}="7-1"
                    !{(12,2)}*+{}*[bricks:c]\frm{**}="7-2"
                    !{(11,2)}*+{}*[bricks:c]\frm{**}="7-3"
                    !{(10,2)}*+{}*[bricks:c]\frm{**}="7-4"
                    !{(9,2)}*+{}*[bricks:c]\frm{**}="7-5"
                    !{(8,2)}*+{}*[bricks:c]\frm{**}="7-6"
                    !{(12,4)}*+{\dots\dots}="dots2"
                    !{(11,4)}*+{\dots\dots}="dots3"
                    !{(10,4)}*+{\dots\dots}="dots4"
                    !{(9,4)}*+{\dots\dots}="dots5"
                    !{(8,4)}*+{\dots\dots}="dots6"
                    !{(7,4)}*+{\dots\dots}="dots7"
                    !{(10,6)}*+{}*[bricks]\frm{**}="4-0"
                    !{(9,6)}*+{}*[bricks:c]\frm{**}="4-1"
                    !{(8,6)}*+{}*[bricks:c]\frm{**}="4-2"
                    !{(7,6)}*+{}*[bricks:c]\frm{**}="4-3"
                    !{(8,8)}*+{}*[bricks:c]\frm{**}="3-0"
                    !{(7,8)}*+{}*[bricks:c]\frm{**}="3-1"
                    !{(6,8)}*+{}*[bricks:c]\frm{**}="3-2"
                    !{(4,12)}*+{}*[bricks]\frm{**}="2-0"
                    !{(3,12)}*+{}*[bricks:c]\frm{**}="2-1"
                    !{(2,14)}*+{}*[bricks:c]\frm{**}="1-0"
                    !{(0,17)}*+{\frac{R}{(\ph_{2h})^\lambda}}="topb"
                    !{(0, 16.3)}*+{}="top_E"
                    !{(0.5, 16.5)}*+{}="top_NE"
                    !{(1, 17)}*+{}="top_N"
                    !{(16,-2)}*+{}="y"
                    !{(-2,18)}*+{}="x"
                    !{(-2,-2.5)}*+{}="o1"
                    !{(-3,-2)}*+{}="o2"
                     !{(-3,0)}*+{G^1_{2 \lambda h}}="Cpari"
                    !{(-3,2)}*+{G^1_{(2 \lambda -1)h +1}}="Cdisp"
                    !{(-3,6)}*+{G^1_{4h}}="C4"
                    !{(-3,8)}*+{G^1_{3h+1}}="C3"
                    !{(-3,12)}*+{G^1_{2h}}="C2"
                    !{(-3,14)}*+{G^1_{h+1}}="C2"
                    "o1":"x"  
                    "o2":"y"  
                     "8-0":@/^2cm/^{d_{2\lambda h}}"top_N"
                     "4-0":@/^1.2cm/^{d_{4h}}"top_NE"
                     "2-0":@/^.2cm/^{d_{2h}}"top_E"
}
$$
\caption{Setup for the higher degree terms in the spectral
sequence for $G^1_n$
 in case a)}\label{table:ehigh:anti:case:a}
\end{center}
\end{table}

For case $b)$, the first column in the spectral sequence relevant
for $\ph_{2h}$ is still even. The differential $d_{h-1}$ may be
identified again as multiplication as in formula
\ref{eq:diff:h-1}, but now it vanishes, since the polynomial is
divisible by $\ph_{2h}$.

\begin{table}[ht]
\begin{center}
\def\objectstyle{\scriptstyle}
\UsePatternFile{macpat.xyp}
    \AliasPattern{bricks:c}{mac14}{xymacpat} 
    \AliasPattern{bricks}{mac06}{xymacpat}
\def\myob {\phantom{o_{a,c}}}
$$
\xygraph{  
                    !{<0cm, 0cm>; <0cm, .25cm>:<.55cm, 0cm>::}
%
                    !{(16,0)}*+{}*[bricks]\frm{**}="8-0"
                    !{(15,0)}*+{}*[bricks]\frm{**}="8-1"
                    !{(14,0)}*+{}*[bricks]\frm{**}="8-2"
                    !{(13,0)}*+{}*[bricks]\frm{**}="8-3"
                    !{(12,0)}*+{}*[bricks]\frm{**}="8-4"
                    !{(11,0)}*+{}*[bricks]\frm{**}="8-5"
                    !{(10,0)}*+{}*[bricks]\frm{**}="8-6"
                    !{(9,0)}*+{}*[bricks]\frm{**}="8-7"
                    !{(14,2)}*+{}*[bricks]\frm{**}="7-0"
                    !{(13,2)}*+{}*[bricks]\frm{**}="7-1"
                    !{(12,2)}*+{}*[bricks]\frm{**}="7-2"
                    !{(11,2)}*+{}*[bricks]\frm{**}="7-3"
                    !{(10,2)}*+{}*[bricks]\frm{**}="7-4"
                    !{(9,2)}*+{}*[bricks]\frm{**}="7-5"
                    !{(8,2)}*+{}*[bricks]\frm{**}="7-6"
                    !{(10,6)}*+{}*[bricks]\frm{**}="6-0"
                    !{(9,6)}*+{}*[bricks]\frm{**}="6-1"
                    !{(8,6)}*+{}*[bricks]\frm{**}="6-2"
                    !{(7,6)}*+{}*[bricks]\frm{**}="6-3"
                    !{(6,6)}*+{}*[bricks]\frm{**}="6-4"
                    !{(5,6)}*+{}*[bricks]\frm{**}="6-5"
                    !{(9,7.5)}*+{\dots\dots}="dots2"
                    !{(8,7.5)}*+{\dots\dots}="dots3"
                    !{(7,7.5)}*+{\dots\dots}="dots4"
                    !{(6,7.5)}*+{\dots\dots}="dots5"
                    !{(5,7.5)}*+{\dots\dots}="dots6"
                    !{(4,7.5)}*+{\dots\dots}="dots7"
                    !{(8,9)}*+{}*[bricks]\frm{**}="3-0"
                    !{(7,9)}*+{}*[bricks]\frm{**}="3-1"
                    !{(6,9)}*+{}*[bricks]\frm{**}="3-2"
                    !{(4,13)}*+{}*[bricks]\frm{**}="2-0"
                    !{(3,13)}*+{}*[bricks]\frm{**}="2-1"
                    !{(2,15)}*+{}*[bricks]\frm{**}="1-0"
                    !{(0,18)}*+{\frac{R}{(\ph_{2h})^{\lambda+1}}}="topb"
                    !{(16,-3)}*+{}="y"
                    !{(-2,18)}*+{}="x"
                    !{(-2,-3.5)}*+{}="o1"
                    !{(-3,-3)}*+{}="o2"
                     !{(-3,0)}*+{G^1_{2 \lambda h}}="C8"
                    !{(-3,2)}*+{G^1_{(2 \lambda -1)h +1}}="C7"
                    !{(-3,6)}*+{G^1_{2 (\lambda-1) h}}="C6"
                    !{(-3,9)}*+{G^1_{3h+1}}="C3"
                    !{(-3,13)}*+{G^1_{2h}}="C2"
                    !{(-3,15)}*+{G^1_{h+1}}="C1"
                    "o1":"x"  
                    "o2":"y"  
                     "7-1":"6-0"
                     "7-2":"6-1"
                     "7-3":"6-2"
                     "7-4":"6-3"
                     "7-5":"6-4"
                     "7-6":_{d_{h+1}}"6-5"
                     "3-1":"2-0"
                     "3-2":_{d_{h+1}}"2-1"
                     "1-0":_{d_{h+1}}"topb"
}
$$
\caption{$E_{h-1}$-term of the spectral sequence for $G^1_n$ in
case b)}\label{table:eh-1:anti:case:b}
\end{center}
\end{table}

The next non-vanishing differential is $d_{h+1}$. See Table
\ref{table:eh-1:anti:case:b}. It takes the module in positive
codimension in an odd column $G^1_{(2 \mu +1)h+1}$ to the elements
in the even column $G^1_{2 \mu h}$ (for $1\leq \mu \leq
\lambda-1$). Via generators, it may be identified with the
multiplication by
\begin{equation}\label{eq:diff:h+1}
\qbin{n-2\mu h}{h+1}=\qbin{2(\lambda-\mu)h+h+1}{h+1}
\end{equation}
and it is therefore injective when restricted to modules in
positive codimension in  odd columns. Further $d_{h+1}$ is also
non-zero as a map $E_{h+1}^{(2 \lambda h-1, h+1)}\to
E_{h+1}^{(n-1,1)}$. Actually the term
$$
E_{h+1}^{(n-1,1)}\cong
E_2^{(n-1,1)}\cong\frac{R}{[2(n-1)]!!/[n-1]!}
$$
has  $R/(\ph_{2h})^{\lambda+1}$ as the only summand with torsion
of type $\ph_{2h}
^l$. It is easy to check that the relative map can be identified
with the multiplication by $\ph_{2h}^\lambda$.

Thus, the only survivors in the $E_{2h}$ term are the first even
column,  the top modules in the odd columns, generated in
positions $(2(\lambda-\mu)h-1 , (2\mu+1)h+1 )$ by $o_{\mu,0}$ for
$1\leq \mu\leq \lambda-1$, as well as $ E_{2h}^{(n-1,1)}$ which
has  $R/(\ph_{2h})^{\lambda}$
 as summand.\\
 Note that the higher
 differentials  vanish when restricted to the first even
 column. Actually we may lift the generators of type $e_{\lambda-s,2s}[h]$ to global generators $e_{\lambda-s,2s+1}[h]$ for $0\leq s < \lambda$.
 Similarly for $0\leq s < \lambda$ we may lift $o_{\lambda-s-1,2s+1}[h]$ to the global generator
 $o_{\lambda-s-1,2s+2}[h]$.
\begin{table}[ht!]
\begin{center}
\def\objectstyle{\scriptstyle}
\UsePatternFile{macpat.xyp}
    \AliasPattern{bricks:c}{mac14}{xymacpat} 
    \AliasPattern{bricks}{mac06}{xymacpat}
\def\myob {\phantom{o_{a,c}}}
$$
\xygraph{  
                    !{<0cm, 0cm>; <0cm, .25cm>:<.55cm, 0cm>::}
%
                    !{(16,0)}*+{}*[bricks]\frm{**}="8-0"
                    !{(15,0)}*+{}*[bricks]\frm{**}="8-1"
                    !{(14,0)}*+{}*[bricks]\frm{**}="8-2"
                    !{(13,0)}*+{}*[bricks]\frm{**}="8-3"
                    !{(12,0)}*+{}*[bricks]\frm{**}="8-4"
                    !{(11,0)}*+{}*[bricks]\frm{**}="8-5"
                    !{(10,0)}*+{}*[bricks]\frm{**}="8-6"
                    !{(9,0)}*+{}*[bricks]\frm{**}="8-7"
                    !{(14,2)}*+{}*[bricks]\frm{**}="7-0"
                    !{(13,2)}*+{}*[bricks:c]\frm{**}="7-1"
                    !{(12,2)}*+{}*[bricks:c]\frm{**}="7-2"
                    !{(11,2)}*+{}*[bricks:c]\frm{**}="7-3"
                    !{(10,2)}*+{}*[bricks:c]\frm{**}="7-4"
                    !{(9,2)}*+{}*[bricks:c]\frm{**}="7-5"
                    !{(8,2)}*+{}*[bricks:c]\frm{**}="7-6"
                    !{(10,6)}*+{}*[bricks:c]\frm{**}="6-0"
                    !{(9,6)}*+{}*[bricks:c]\frm{**}="6-1"
                    !{(8,6)}*+{}*[bricks:c]\frm{**}="6-2"
                    !{(7,6)}*+{}*[bricks:c]\frm{**}="6-3"
                    !{(6,6)}*+{}*[bricks:c]\frm{**}="6-4"
                    !{(5,6)}*+{}*[bricks:c]\frm{**}="6-5"
                    !{(9,7.5)}*+{\dots\dots}="dots2"
                    !{(8,7.5)}*+{\dots\dots}="dots3"
                    !{(7,7.5)}*+{\dots\dots}="dots4"
                    !{(6,7.5)}*+{\dots\dots}="dots5"
                    !{(5,7.5)}*+{\dots\dots}="dots6"
                    !{(4,7.5)}*+{\dots\dots}="dots7"
                    !{(8,9)}*+{}*[bricks]\frm{**}="3-0"
                    !{(7,9)}*+{}*[bricks:c]\frm{**}="3-1"
                    !{(6,9)}*+{}*[bricks:c]\frm{**}="3-2"
                    !{(4,13)}*+{}*[bricks:c]\frm{**}="2-0"
                    !{(3,13)}*+{}*[bricks:c]\frm{**}="2-1"
                    !{(2,15)}*+{}*[bricks:c]\frm{**}="1-0"
                    !{(0,18)}*+{\frac{R}{(\ph_{2h})^{\lambda}}}="topb"
                     !{(0.5, 17.5)}*+{}="top_NE"
                    !{(1, 18)}*+{}="top_N"
                    !{(16,-3)}*+{}="y"
                    !{(-2,18)}*+{}="x"
                    !{(-2,-3.5)}*+{}="o1"
                    !{(-3,-3)}*+{}="o2"
                     !{(-3,0)}*+{G^1_{2 \lambda h}}="C8"
                    !{(-3,2)}*+{G^1_{(2 \lambda -1)h +1}}="C7"
                    !{(-3,6)}*+{G^1_{2 (\lambda-1) h}}="C6"
                    !{(-3,9)}*+{G^1_{3h+1}}="C3"
                    !{(-3,13)}*+{G^1_{2h}}="C2"
                    !{(-3,15)}*+{G^1_{h+1}}="C1"
                    "o1":"x"  
                    "o2":"y"  
                    "7-0":@/^2cm/^{d_{(2 \lambda-1)h +1}}"top_N"
                    "3-0":@/^1cm/^{d_{3h +1}}"top_NE"
}
$$
\caption{Setup for the higher degree terms in the spectral
sequence for $G^1_n$ in case b)}\label{table:higher:anti:case:b}
\end{center}
\end{table}
Finally, as in case a), the module in positions
$(2(\lambda-\mu)h-1 , (2\mu+1)h+1 )$ for
$1\leq \mu\leq \lambda -1$ vanish in the higher terms of the
spectral sequence while the module in position $(n-1,1)$ has
eventually as summand $R/\ph_{2h}$. Clearly the coboundary
$o_{\lambda,0}[h]$ projects onto a generator of the latter.

Case c) and d) present no new complications and are omitted. \qed
\end{dm}


\subsection{Spectral sequence for $\br{\t{B}_n}$}
We can now compute the cohomology $H^*(G_{\widetilde{B}_n}, R_{q,t})$.
We will do this by means of the Salvetti complex $\widehat{C}^*\widetilde{B}_n$.

As in Section [\ref{s:splitting}], let $\widehat{\inv}{\widetilde{B}_n}$ be the module of
the $\sigma$-invariant elements and $\widehat{\anti}{\widetilde{B}_n}$ the module of the $\sigma$-anti-invariant elements.
We can split our module $\widehat{C^*}{\widetilde{B}_n}$ into the direct sum:
$$
\widehat{C^*}{\widetilde{B}_n} = \widehat{\inv}{\widetilde{B}_n} \oplus \widehat{\anti}{\widetilde{B}_n}.
$$
Using the map $\beta: C^*{B_n} \to \widehat{C^*}{\widetilde{B}_n}$ so defined:
$$\beta: 0A \mapsto \abin{0}{0}A$$
$$\beta: 1A \mapsto \abin{1}{0}A + \abin{0}{1}A$$
one can see that the submodule $\widehat{\inv}{\widetilde{B}_n}$
is isomorphic (as a differential complex) to $C^*{B_n}$. Its
cohomology has been computed in \cite{cohom_B_n}. We recall the
result:
\begin{theorem} [\cite{cohom_B_n}] \label{t:inv}
$$
H^i(G_{B_n}, R_{q,t})= \left\{
\begin{array}{ll}
\bigoplus_{d \mid n, 0 \leq i \leq d-2}
\{ d \}_i \oplus \{ 1 \}_{n-1} & \mbox{ if } i = n\\
\bigoplus_{
d \mid n, 0 \leq i \leq d-2,
d \leq \frac{n}{j+1}
}
\{ d \}_i & \mbox{ if }i = n -2j \\
\bigoplus_{d \nmid n,
d \leq \frac{n}{j+1} }
\{ d \}_{n-1} & \mbox{ if }i = n -2j -1.
\end{array}
\right.
$$ \qed
\end{theorem}
Hence we only need to compute the cohomology of
$\widehat{\anti}{\widetilde{B}_n}$. In order to do this we make
use of the results presented in Section \ref{ss:anti-inv}.
First consider the subcomplex of $\widehat{C^*}{\widetilde{B}_n}$
defined as
$$
\oC_n = <\abin{0}{1}A, \abin{1}{1}A>.
$$
We define the map $\kappa: \oC_n \to \widehat{\anti}{\widetilde{B}_n}$ by
$$
\kappa: \abin{0}{1}A \mapsto \abin{0}{1}A -  \abin{1}{0}A
$$
$$
\kappa: \abin{1}{1}A \mapsto 2 \abin{1}{1}
A.
$$
It is easy to check that $\kappa$ gives an isomorphism of differential complex. Now we define a filtration $\cF$ on the complex $\oC_n$:

$$
\cF_i \oC_n = < \abin{0}{1}A1^i, \abin{1}{1}A1^i>.
$$

The quotient $\cF_i \oC_n / \cF_{i+1} \oC_n$ is isomorphic to the
complex $\left(G^1_{n-i}[t^\pmu]\right)
[i]$ (see Proposition \ref{p:dpssG}) with trivial action on the
variable $t$. Hence we use the spectral sequence defined by the
filtration $\cF$ to compute the cohomology of the complex $\oC_n$.

The $E_0$-term of the spectral sequence is given by
\begin{align*}
E_0^{i,j}
=&  \frac{\left(\cF_i \oC_n\right)^{(i+j)} }{ \left(\cF_{i+1}
\oC_n\right)^{(i+j)}}\\
=&\left((G^1_{n-i})^{(i+j)}[t^\pmu]\right)[i]\\
%
%
%
=&  (G^1_{n-i})^j   [t^\pmu]
\end{align*}
for $0 \leq i \leq n-2$. Finally:
\begin{align*}
E_0^{n-1,1} &= R &   E_0^{n,1} &= R
\end{align*}and all the other terms are zero. The differential $d_0: E_0^{i,j}
\to 
E_0^{i,j+1}$ corresponds to the differential on the complex
%
%
$G^1_{n-i}$.It follows that the $E^1$-term is given by the
cohomology of the complexes
$G^1_{n-i}$:
$$
E^1_{i,j} = H^j(G^1_{n-i})[t^\pmu]
$$
for $0 \leq i \leq n-2$ and
$$
E_1^{n-1,1} = R, \quad E_1^{n,1} = R.
$$


As in Section \ref{ss:anti-inv}, we can separately consider in the spectral sequence
$E_*$ the modules with torsion of type $\ph_{2h}^l$ for an integer $h \geq 1$.

For a fixed integer $h >0$, let $c \in \{ 0, \ldots, 2h-1 \}$ be the congruency class of $n$ $\mod{2h}$
and let $\lambda$ be an integer such that $n = c + 2 \lambda h$. We consider the two cases:

a) $0 \leq c \leq h$;

b) $h+1 \leq c \leq 2h-1$.

In case a) the modules of $\ph_{2h}$-torsion are:\\
with $0 \leq \mu \leq \lambda -1, 0 \leq i \leq \lambda - \mu -1$
$$
E_1^{c +2 \mu h, 2 (\lambda - \mu)h - 2i} \simeq \{ 2h \}[t^\pmu]
$$
generated by $e_{\lambda - \mu-i,2i}[h]01^{c +2 \mu h}$;\\
with $0 \leq \mu \leq \lambda -1, 0 \leq i \leq \lambda - \mu -1 $
$$
E_1^{c +2 \mu h, 2 (\lambda - \mu)h - 2i - 1} \simeq \{ 2h \}[t^\pmu]
$$
generated by $o_{\lambda - \mu -i-1,2i+1}[h]01^{c +2 \mu h}$; \\
with $0 \leq \mu \leq \lambda-1, 0 \leq i \leq \lambda - \mu-1 $
$$
E_1^{c + 2 \mu h + h -1, 2 (\lambda - \mu)h -h +1 - 2i } \simeq \{ 2h \}[t^\pmu]
$$
generated by $o_{\lambda - \mu -i-1,2i}[h]01^{c + 2 \mu h + h-1} $; \\
with $0 \leq \mu \leq \lambda-2, 0 \leq i \leq \lambda - \mu -2 $
$$
E_1^{c + 2 \mu h + h -1, 2 (\lambda - \mu)h -h +1 - 2i -1 } \simeq \{ 2h \}[t^\pmu]
$$
generated by $e_{\lambda - \mu -i-1,2i+1}[h]01^{c + 2 \mu h + h-1} $.

In case b) the modules of $\ph_{2h}$-torsion are: \\ 
with $0 \leq
\mu \leq \lambda -1, 0 \leq i \leq \lambda - \mu -1$
$$
E_1^{c +2 \mu h, 2 (\lambda - \mu)h - 2i} \simeq \{ 2h \}[t^\pmu]
$$
     generated by $e_{\lambda - \mu-i  ,2i}[h]01^{c +2 \mu h }$; \\
%
%
%
with $0 \leq \mu \leq \lambda -1, 0 \leq i \leq \lambda - \mu -1$
$$
E_1^{c +2 \mu h, 2 (\lambda - \mu)h - 2i - 1} \simeq \{ 2h \}[t^\pmu]
$$
       generated by $o_{\lambda - \mu -i -1 ,2i+1}[h]01^{c +2 \mu h}$; \\
%
%
%
%
%
with $0 \leq \mu \leq \lambda, 0 \leq i \leq \lambda - \mu$
$$
E_1^{c + 2 \mu h - h -1, 2 (\lambda - \mu)h +h +1 - 2i } \simeq \{ 2h \}[t^\pmu]
$$
generated by $o_{\lambda - \mu -i,2i}[h]01^{c + 2 \mu h - h-1} $;\\
with $0 \leq \mu \leq \lambda-1, 0 \leq i \leq \lambda - \mu -1$
$$
E_1^{c + 2 \mu h - h -1, 2 (\lambda - \mu)h +h +1 - 2i -1 } \simeq \{ 2h \}[t^\pmu]
$$
generated by $e_{\lambda - \mu -i,2i+1}[h]01^{c + 2 \mu h - h-1} $.

In the $E_1$-term of the spectral sequence, the only non-trivial map is the map $d_1:E_1^{n-1,1} \to E_1^{n,1}$, that corresponds to the multiplication by the polynomial
$$
\frac{{\widehat{W}}_{\widetilde{B}_n}[q,t]}
{W_{B_n}[q,t]} = \prod_{i=1}^{n-1}(1+q^{i}) = \prod_{h \leq n} \ph_{2h}^{\lfloor\frac{n-1}{h}\rfloor-\lfloor\frac{n-1}{2h}\rfloor}.
$$
Then in $E_2$ we have:
$$
E_2^{n-1,1} = 0
$$
and
$$
E_2^{n,1} = \bigoplus R/(\ph_{2h}^{\lfloor\frac{n-1}{h}\rfloor-\lfloor\frac{n-1}{2h}\rfloor}).
$$\\
Notice that the integer $f(n,h) = \lfloor\frac{n-1}{h}\rfloor-\lfloor\frac{n-1}{2h}\rfloor$ corresponds to $\lambda$ in case a) and to $\lambda +1$ in case b).

Now we consider the higher
differentials in the spectral sequence. The first possibly
non-trivial maps are $d_{h-1}$ and $d_{h+1}$. In case a) the map
$d_{h-1}$ is given by the multiplication by
$$
\prod_{i=n}^{n+h-2}(1+tq^i)
$$
and the map $d_{h+1}$ is the null map.
The maps
$$
d_{2(\lambda - \mu)h}:\{2h \}[t^\pmu]= E_{2(\lambda - \mu)h}^{c+2\mu h, 2(\lambda-\mu)h} \to E_{2(\lambda - \mu)h}^{n,1}
$$
where $\mu$ goes from $\lambda-1$ to $0$, correspond, up to invertibles, modulo $\ph_{2h}$,to multiplication by 
$$
\ph_{2h}^\mu (\prod_{i=0}^{2h-1} (1+tq^i))^{\lambda - \mu}.
$$
Moreover they are all injective and the term $E_{2(\lambda)h+1}^{n,1}$ is given by the 
quotient
$$
R/(\ph_{2h}^\lambda, \ph_{2h}^{\lambda-1} \prod_{i=0}^{2h-1} (1+tq^i), \ldots, (\prod_{i=0}^{2h-1} (1+tq^i))^{\lambda}) =
$$
$$
= R/(\ph_{2h},\prod_{i=0}^{2h-1} (1+tq^i) )^\lambda.
$$

In case b) the map $d_{h-1}$ is null and the map
$d_{h+1}$ is the multiplication by the polynomial
$$
\prod_{i=n+h-1}^{n+2h-1}(1+tq^i).
$$
The maps
$$
d_{2(\lambda - \mu)h + h +1}: \{2h \}[t^\pmu] = E^{c + 2\mu h + h -1, 2(\lambda - \mu)h -h }_{2(\lambda - \mu)h + h +1} \to E_{2(\lambda - \mu)h + h +1}^{1,n}
$$
%
%
where $\mu$ goes from $\lambda $ to $0 $, correspond, up to invertibles, modulo $\ph_{2h}$,to multiplication by 
$$
\ph_{2h}^\mu (\prod_{i=0}^{2h-1} (1+tq^i))^{\lambda - \mu +1}.
$$
Hence they are all injective and the term $E_{2(\lambda)h+h+2}^{n,1}$ is given by the 
quotient
$$
R/(\ph_{2h},\prod_{i=0}^{2h-1} (1+tq^i) )^{\lambda+1}.
$$
Since all the generators lift to global cocycles, it turns out that
all the other differentials are null.
%
%
%
%
%
Hence we proved the following:
\begin{theorem} \label{t:anti}
$$
H^{n+1} (\widehat{\anti}{\widetilde{B}_n}) \simeq \displaystyle \bigoplus_{h>0} \{ \{ 2h \} \}_{f(n,h)} 
$$
and, for $s \geq 0$:

$$
H^{n-s} (\widehat{\anti}{\widetilde{B}_n}) \simeq \bigoplus_{\begin{array}{c} \scriptstyle h>2 \\
\scriptstyle i \in I(n,h)
\end{array}} \{2h\}_i^{\oplus max(0,\lfloor \frac{n}{2h}\rfloor-s)}
$$
with $I(n,h) = \{n, \ldots, n+h-2 \}$ if $n \simeq 0, 1, \ldots, h \mod{2h}$, $ f(n,h)= \lfloor \frac{n+h-1}{2h} \rfloor$ and $I(n,h) = \{n+h-1, \ldots, n+2h-1 \} $ if $n \simeq h+1, h+2, \ldots, 2h-1 \mod{2h}$. \qed
\end{theorem}

Putting together the results of Theorem \ref{t:inv} and \ref{t:anti},
we get Theorem \ref{thm:main2}.

As a corollary, we  use the long exact
sequences associated to
$$
0 \too \Q[[t^\pmu]] \stackrel{m(q)}{\too} M \stackrel{1+q}{\too} M \too 0
$$
and
$$
0 \too \Q \stackrel{m(t)}{\too} \Q[[t^\pmu]] \stackrel{1+t}{\too} \Q[[t^\pmu]] \too 0
$$
to get the constant coefficients cohomology for $G_{\widetilde{B}_n}$.
Here $m(x)$ is the multiplication by the series
$$\sum_{i\in \Z}\ (-x)^i.$$
We give only the result, omitting details which come from non difficult
analysis of the above mentioned sequences and recalling that the Euler characteristic of the complex is $1$, for $n$ even, and $-1$, for $n$ odd.

\begin{theorem}\label{thm:main3}
$$
H^{i}({G}_{\widetilde{B}_n}, \Q)= \left\{
\begin{array}{lcl}
\Q & \mbox{ if } & i = 0 \\
\Q^2 & \mbox{ if } & 1 \leq i \leq n-2 \\
\Q^{2+\lfloor \frac{n}{2}\rfloor} & \mbox{ if } & i = n-1,n 
\end{array}
\right.
$$
where the $t$ and $q$
actions correspond to the multiplication by $-1$. \qed
\end{theorem}

\providecommand{\bysame}{\leavevmode\hbox
to3em{\hrulefill}\thinspace}
\providecommand{\MR}{\relax\ifhmode\unskip\space\fi MR }
\providecommand{\MRhref}[2]{%
  \href{http://www.ams.org/mathscinet-getitem?mr=#1}{#2}
} \providecommand{\href}[2]{#2}

\end{document}